\documentclass[conference]{IEEEtran}
\usepackage[numbers,sort&compress]{natbib}
\usepackage{amsmath,amsthm,amsfonts,amssymb,bm}
\usepackage{subfigure}
\usepackage{booktabs}
\usepackage{geometry}
\usepackage{graphicx}
\usepackage{algorithm}
\usepackage{algorithmicx}
\usepackage{algpseudocode}
\usepackage{multirow}
\usepackage{hyperref}
\usepackage[figuresright]{rotating}

\ifCLASSINFOpdf
\else
\fi
\DeclareRobustCommand*{\IEEEauthorrefmark}[1]{%
	\raisebox{0pt}[0pt][0pt]{\textsuperscript{\footnotesize\ensuremath{#1}}}}

\begin{document}
%
\title{Calico Salmon Migration Algorithm: A novel meta-heuristic optimization algorithm}

\author{\IEEEauthorblockN{Chao Min\IEEEauthorrefmark{1,2,3},
		Junyi Cui\IEEEauthorrefmark{1,2},
		Liwen Zhou\IEEEauthorrefmark{1,2*}, 
		Qian Yin\IEEEauthorrefmark{1},
		Yijia Wang\IEEEauthorrefmark{1,2}}
	\IEEEauthorblockA{\IEEEauthorrefmark{1}School of Sciences,
		Southwest Petroleum University, Sichuan, China}
	\IEEEauthorblockA{\IEEEauthorrefmark{2}Institute for Artificial Intelligence, Southwest Petroleum University, Sichuan, China}
	\IEEEauthorblockA{\IEEEauthorrefmark{3}National Key Laboratory of Oil and Gas Reservoir Geology and exploitation, \\Southwest Petroleum University, Sichuan, China}
	Corresponding Author: Liwen Zhou \quad Email: zhouliwen@live.cn}


%


\maketitle

\begin{abstract}
A novel population-based optimization method is proposed in this paper, the Calico Salmon Migration Algorithm (CSMA), which is inspired by the natural behavior of calico salmon during their migration for mating. The CSMA optimization process comprises four stages: selecting the search space by swimming into the river, expanding the search space from the river into the ocean, performing precise search during the migrating process, and breeding new subspecies by the remaining calico salmon population. To evaluate the effectiveness of the new optimizer, we conducted a series of experiments using different optimization problems and compared the results with various optimization algorithms in the literature. The numerical experimental results for benchmark functions demonstrate that the proposed CSMA outperforms other competing optimization algorithms in terms of convergence speed, accuracy, and stability. Furthermore, the Friedman ranking test shows that the CSMA is ranked first among similar algorithms. 
\end{abstract}

%
\IEEEpeerreviewmaketitle

\section{Introduction}
In this paper, we proposed a novel nature-inspired meta-heuristic optimization algorithm called the Calico Salmon Migration Algorithm (CSMA), which emulates the natural behavior of calico salmon as they migrate to the sea for spawning. The optimization process of CSMA follows the natural spiral evolution process, which allows the best genes within sub-populations to be preserved and enables them to tackle environmental challenges. To verify the robustness and effectiveness of CSMA, we conducted extensive experiments using 23 CEC-2005 test functions. Our results demonstrate that CSMA consistently outperforms other competing algorithms.

\section{Calico Salmon Migration Algorithm (CSMA)}
The details of the proposed natural meta-heuristic algorithm, Calico Salmon Migration Algorithm (CSMA), are described in the subsequent sections.

\subsection{Behavior and Inspiration during Calico Salmon Migration}
North Pacific calico salmon inhabit primarily Russia, Upper Alaska, Japan, and Mexico, playing a crucial role in maintaining the structure, function, and processes of riverine and marine ecosystems \cite{bib1,bib2}. These salmon undergo migration, as they mature in the ocean and then return to freshwater for reproduction, continuing this cycle for generations. From the perspective of bio-mimicry, the competitive survival process of calico salmon migration and spawning can be seen as an analogy to solving an optimization problem. Each new generation of sub-species obtained during this process corresponds to an optimal solution for the given optimization problem. This natural behavior of problem-solving has been abstracted into a mathematical paradigm, which characterizes meta-heuristic algorithms. In addition, researchers have also observed this phenomenon and proposed Salmon Migration Optimization (SMO) \cite{bib3}. However, this article specifically focuses on Water Flow, Magnetic, and Pheromone Oriented Heuristic, which fundamentally differs from the migration discussed in this paper.

In our simulation, we noticed that the exploration phase of the solution space corresponds to the selection and food search behaviors of calico salmon. Conversely, the exploitation phase of the space corresponds to the migration and mating behaviors of calico salmon. The following subsections provide a detailed explanation of how these natural behavioral processes of calico salmon have been integrated into the CSMA algorithm.

\subsubsection{Initialization of the solution space}
CSMA explores the search space similarly to that the baby salmons migrate from the river to the ocean. The first generation of CSMA is a set of uniformly distributed random candidate solutions X, which can be represented as formula (1):
\begin{equation}
	X = \left[ {\begin{array}{*{20}{c}}
			{{x_{1,1}}} &  \cdots  & {{x_{1,j}}} &  \cdots  & {{x_{1,D}}}  \\
			{{x_{2,1}}} &  \cdots  & {{x_{2,j}}} &  \cdots  & {{x_{2,D}}}  \\
			\vdots  &  \cdots  &  \vdots  &  \vdots  &  \vdots   \\
			{{x_{N - 1,1}}} &  \cdots  & {{x_{N - 1,j}}} &  \cdots  & {{x_{N - 1,D}}}  \\
			{{x_{N,1}}} &  \cdots  & {{x_{N,j}}} &  \cdots  & {{x_{N,D}}}  \\
	\end{array}} \right],
\end{equation}
here ${x_{i,j}}$ denotes the ${j^{th}}$ position of the ${i^{th}}$ solution, generated randomly by using formula (2). $D$ is the dimension of the problem, and $N$ is the size of the set of candidate solutions.
\begin{equation}
	\begin{array}{l}
		{x_{i,j}} = \left( {UB - LB} \right) \times rand + LB,\: \\ 
		\;\;\;\;\;\;i = 1,2, \cdots N,\:j = 1,2, \cdots D, \\ 
	\end{array}
\end{equation}
where $ rand $ is a random number between 0 and 1, $ UB $ and $ LB $ denote the upper and lower bounds of the given problem.

\subsection{Mathematical model of CSMA}

The proposed CSMA simulates the various stages of the calico salmon life cycle, from fertilized eggs to juveniles and eventually to adults for spawning. Consequently, the optimization process of CSMA is divided into four stages:

(1) $ Selecting $ stage, in which randomly selected candidate solutions are initially searched in the search space; 

(2) $ Expanding\;search $ stage, in which the algorithm further expands the search space;

(3) $ Migrating $ stage, in which the algorithm narrows the search space for a more efficient search based on the exploration phase; 

(4) $ Mating $ stage, in which the algorithm increases the exploitation and searches for the optimal solution/near-optimal solution.

We model the behavior of calico salmon as a mathematical optimization paradigm that identifies the best solution under a given constraint.

\subsubsection{Population energy}

In this paper, we propose a general mathematical rule of the change of the calico salmon population size, which is later called $ Population\;Energy $. It can be defined as a segmented function, as shown in Eq. (3).
\begin{equation}
Ene = \left\{ \begin{array}{l}
	{E_1} = {\left( {\frac{{2{t^2}}}{{{T^2}}}} \right)^{\left( {1 - \frac{t}{T}} \right)}} \\ 
	t \le \frac{2}{3}T; \\ 
	{E_2} = {\left( {1 - \frac{{3\left( {t - \frac{7}{{10}}T} \right)}}{T}} \right)^{\left( {\frac{{30\left( {t - \frac{7}{{10}}T} \right)}}{T}} \right)}} \\ 
	otherwise. \\ 
\end{array} \right.
\end{equation}

Population energy is an adaptive parameter in CSMA that is vital for maintaining a balance between exploration and exploitation in each iteration of the algorithm. It encompasses two aspects: $ E_{1} $ and $ E_{2} $. During the exploration stage, $ E_{1} $ regulates the parameters and is primarily responsible for exploring the search space, identifying new regions, and potential solutions. On the other hand, during the exploitation stage, $ E_{2} $ regulates the parameters and focuses on developing promising areas and improving the solutions to enhance overall performance.

\subsubsection{Selecting stage}
The first stage of calico salmon development involves their transformation from fertilized eggs to juveniles, during which their main focus is on foraging for food. This behavior can be mathematically expressed by Eq. (4).
\begin{equation}
	\begin{array}{c}
		{x_1}\left( {t + 1} \right) = {x_{best}}\left( t \right) \times {E_1} \\ 
		+ \left( {{x_{best}}\left( t \right) - {x_R}\left( t \right) \times {R_B}} \right) \times rand, \\ 
	\end{array}	
\end{equation}
where $ x_{1}(t+1) $ represents the solution generated in the next iteration $ t+1 $ by the first stage $ (x_{1}) $. $ E_{1} $ represents a small value that primarily regulates the search process of the calico salmon in the early stages. $ x_{best}(t) $ denotes the optimal solution obtained before the $ t^{th} $ iteration, which reflects the approximate position of the calico salmon. $ x_{R}(t) $ represents a random solution selected within the range of $ \left[ {1,N} \right] $ at the ${t^{th}}$ iteration. The standard Brownian motion $ (R_{B}) $ is a stochastic process where the step length is sampled from a probability function defined by the Normal (Gaussian) distribution with zero mean $ \mu=0 $ and unit variance $ \sigma ^2 =1 $. The Probability Density Function (PDF) governing this motion at point $ x $ is expressed as Eq. (5) \cite{bib4}.
\begin{equation}
	R_{B}\left( {x;\mu ,\sigma } \right) = \frac{1}
	{{\sqrt {2\pi } }}\exp \left( { - \frac{{{x^2}}}
		{2}} \right).
\end{equation}

\subsubsection{Expanding search stage}
\begin{figure*}[h]
	\centering
	\subfigure[2D-Brownian motion]{
		\begin{minipage}[t]{0.31\linewidth}
			\centering
			\includegraphics[width=1.1\textwidth]{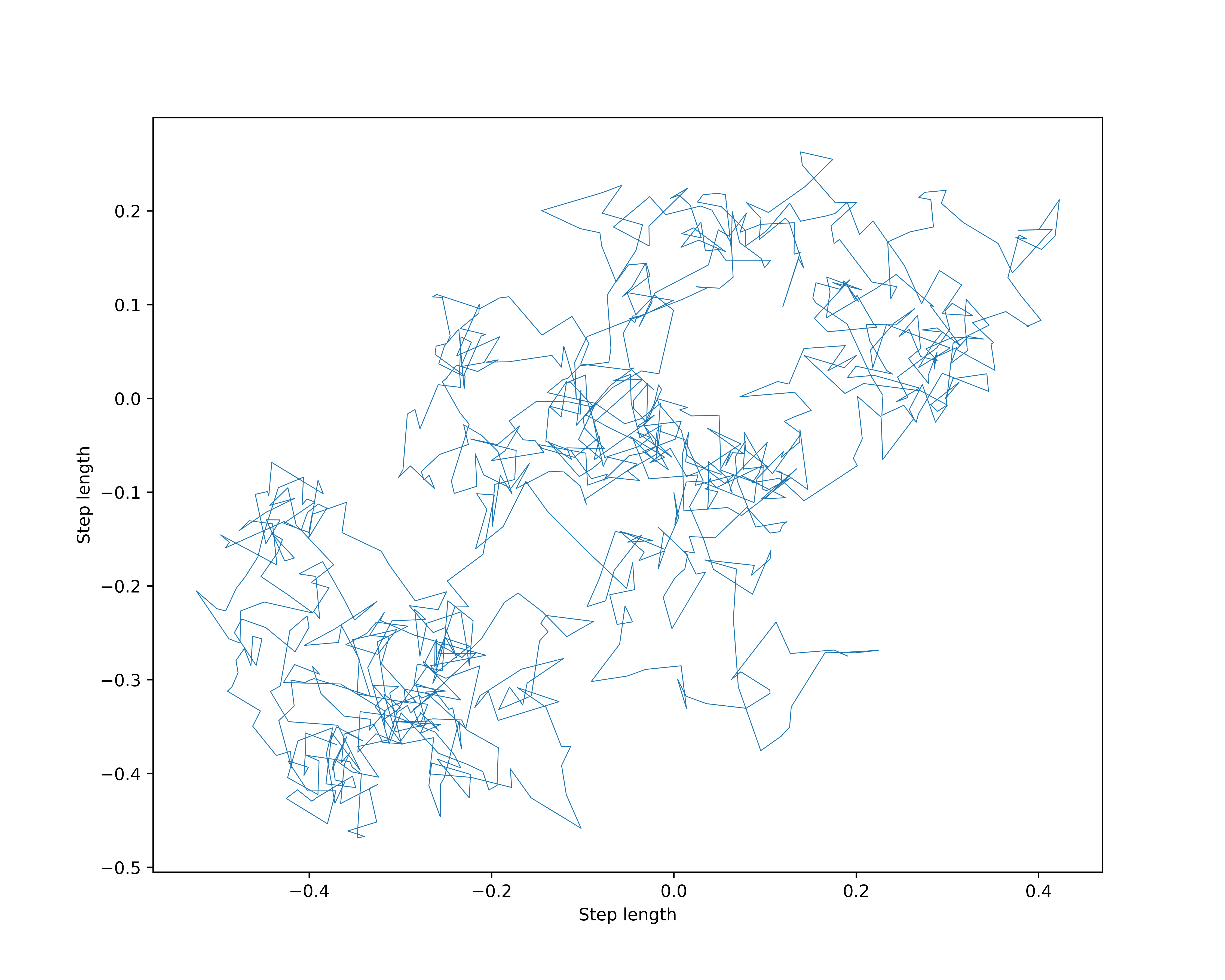}
		\end{minipage}
	}%
	\subfigure[3D-Brownian motion]{
		\begin{minipage}[t]{0.31\linewidth}
			\centering
			\includegraphics[width=1.1\textwidth]{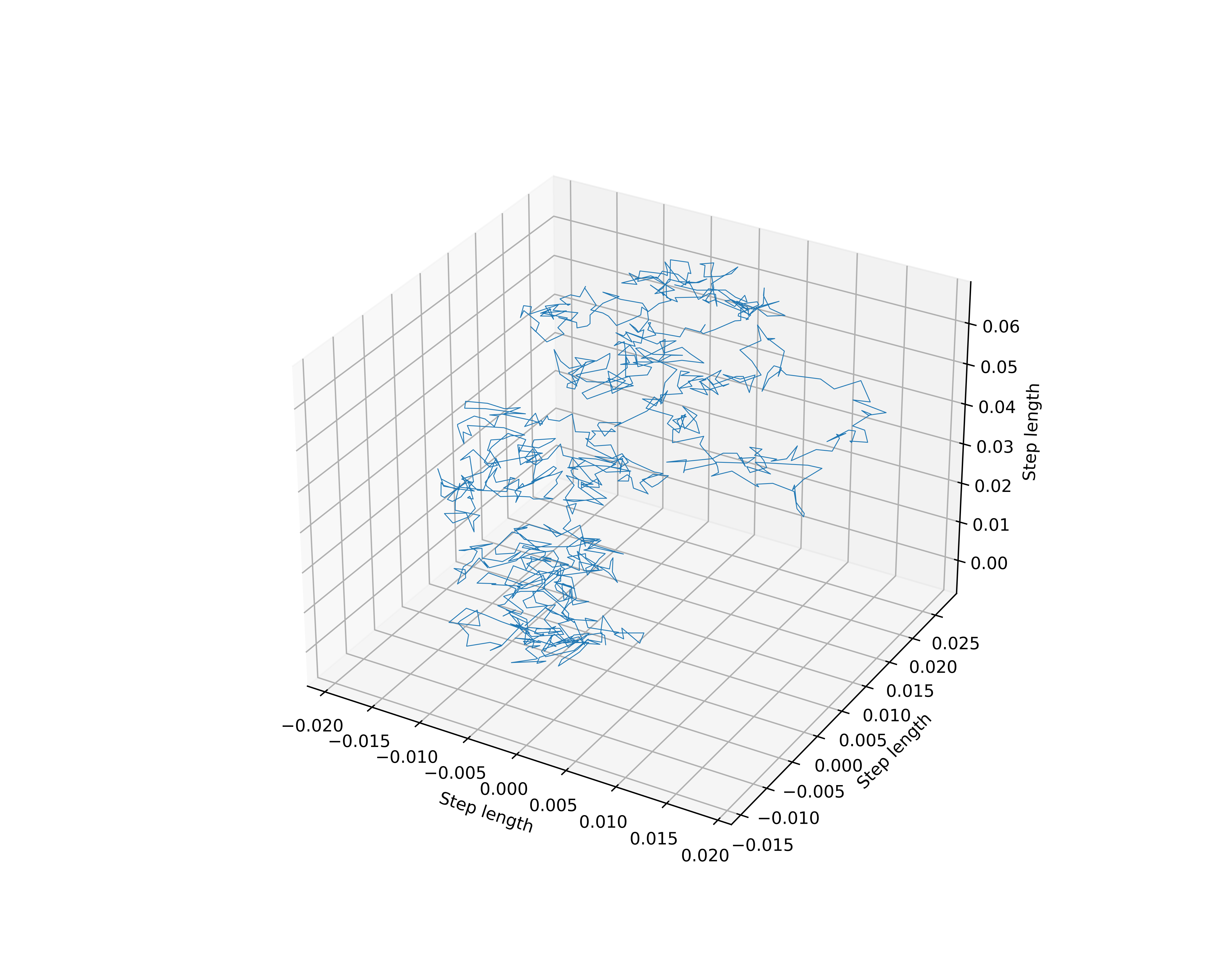}
		\end{minipage}
	}%
	\subfigure[Brownian distribution]{
		\begin{minipage}[t]{0.31\linewidth}
			\centering
			\includegraphics[width=1.1\textwidth]{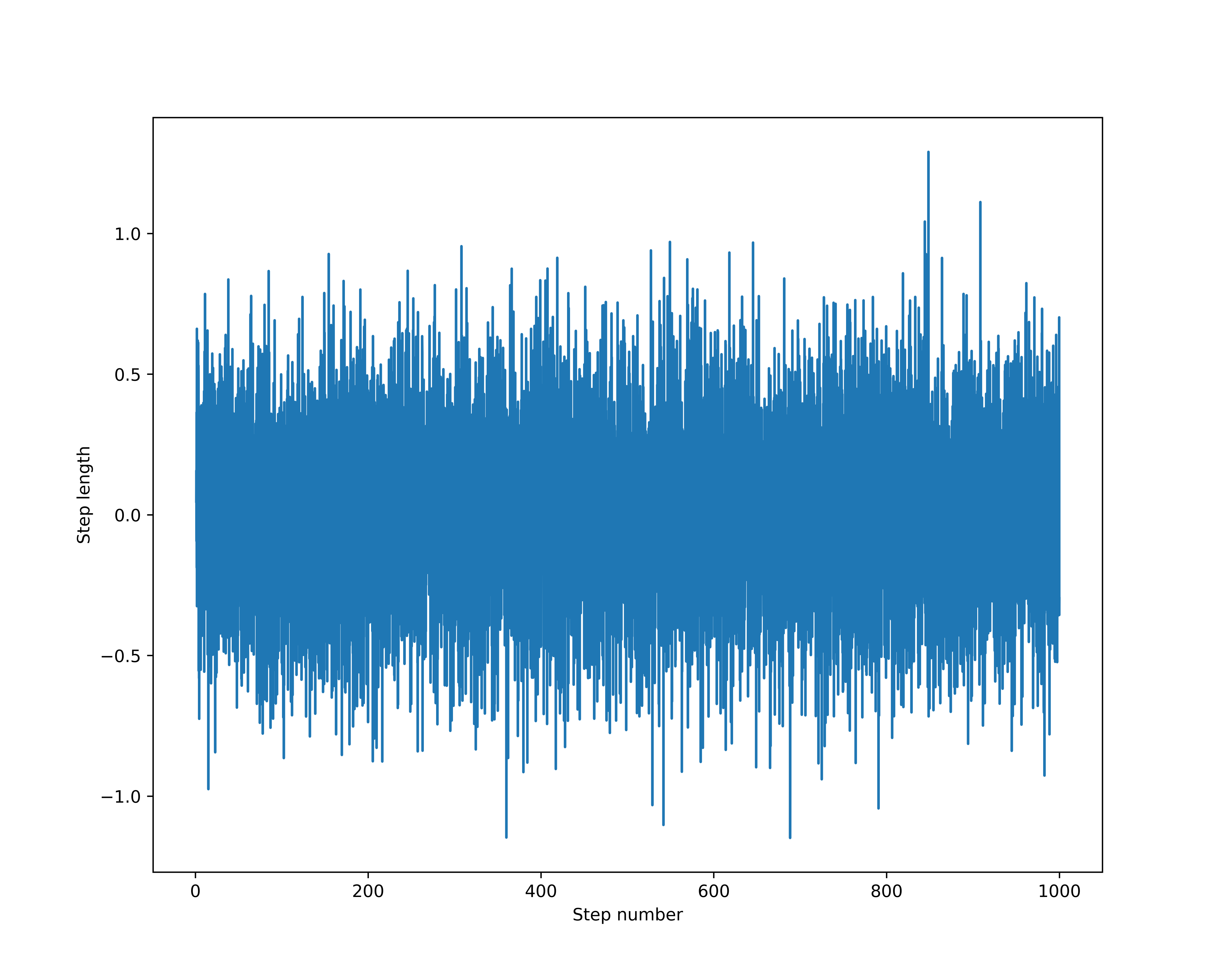}
		\end{minipage}
	}%
	
	\subfigure[2D-L\'{e}vy Flight]{
		\begin{minipage}[t]{0.31\linewidth}
			\centering
			\includegraphics[width=1.1\textwidth]{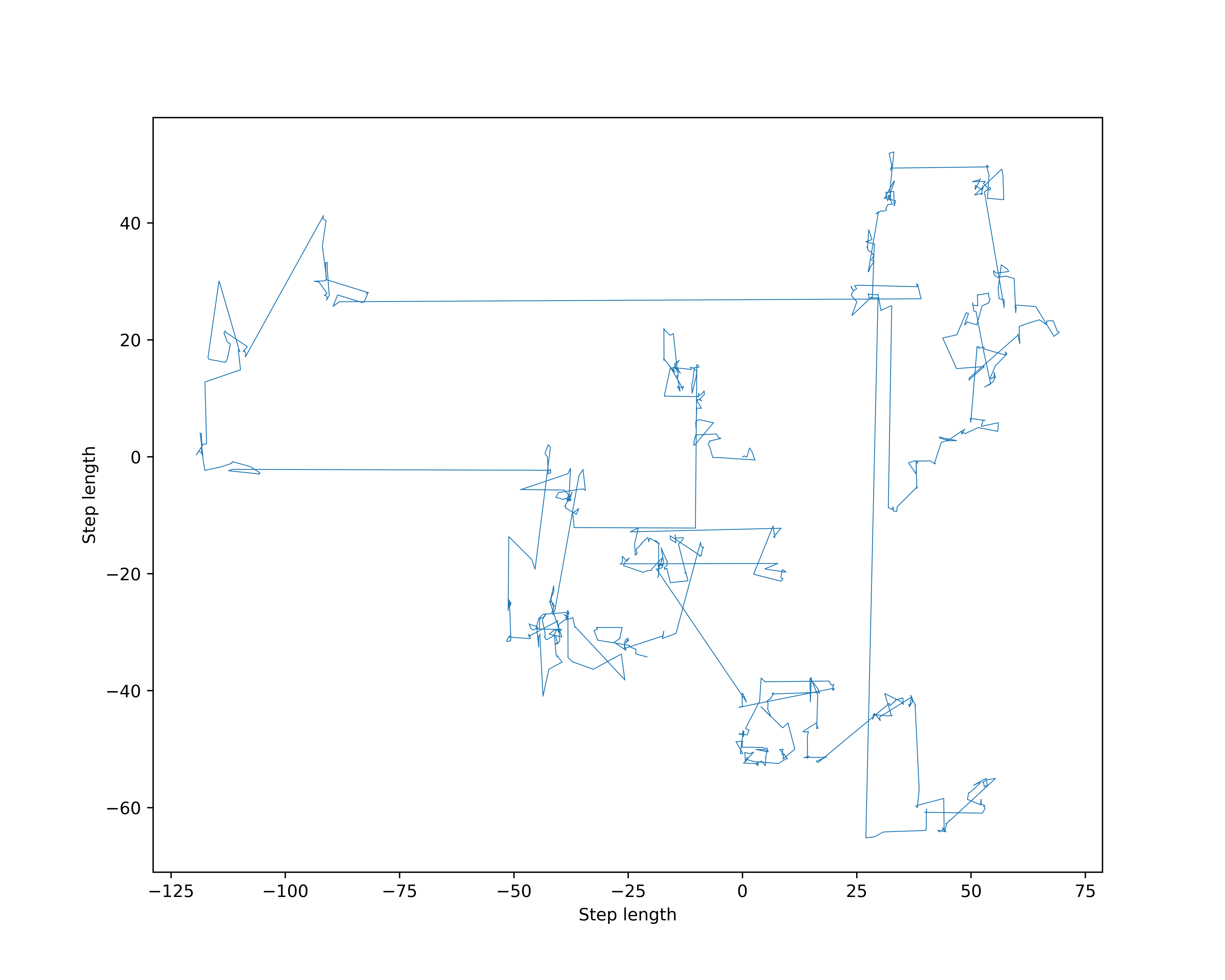}
		\end{minipage}
	}%
	\subfigure[3D-L\'{e}vy Flight]{
		\begin{minipage}[t]{0.31\linewidth}
			\centering
			\includegraphics[width=1.1\textwidth]{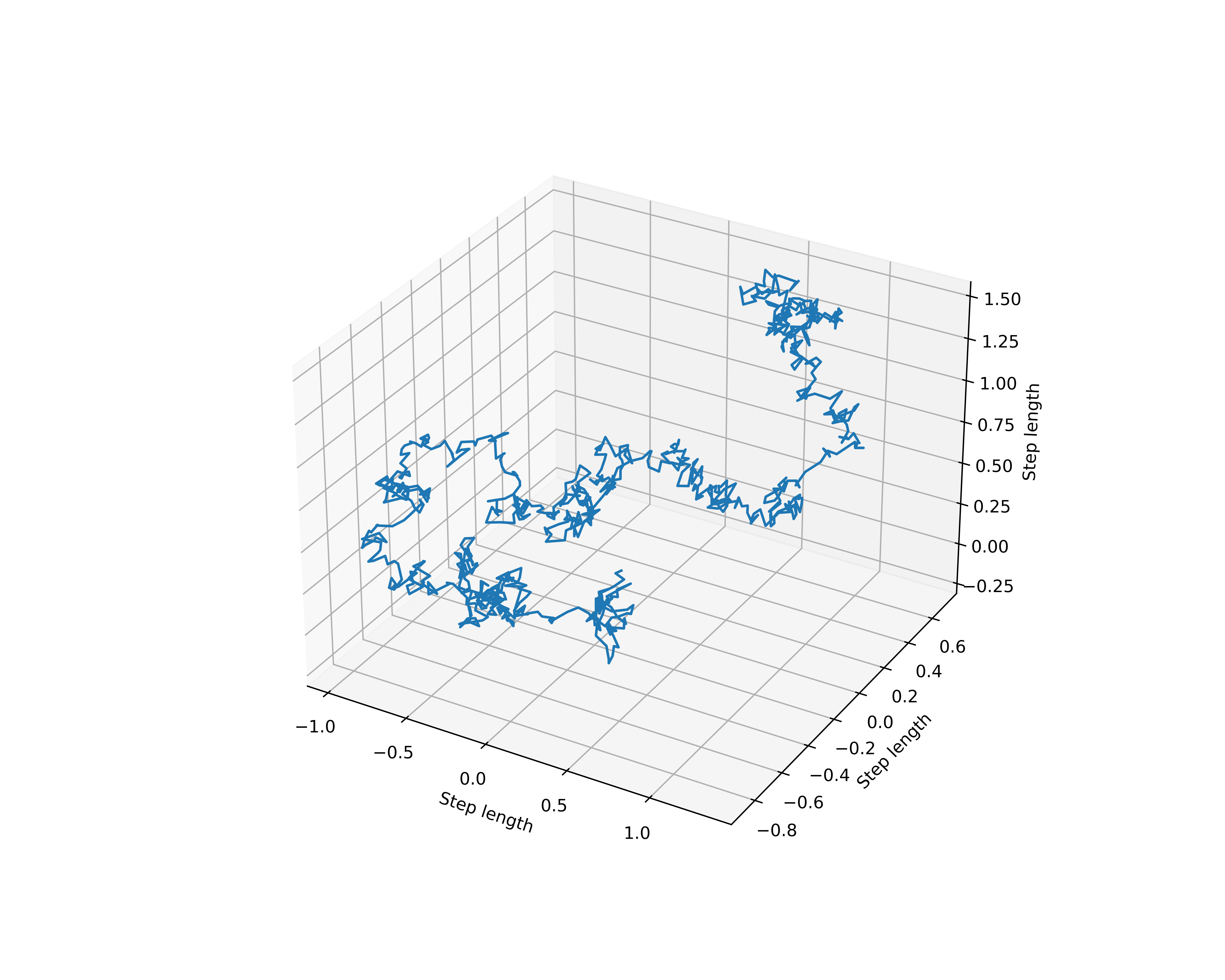}
		\end{minipage}
	}%
	\subfigure[L\'{e}vy distribution]{
		\begin{minipage}[t]{0.31\linewidth}
			\centering
			\includegraphics[width=1.1\textwidth]{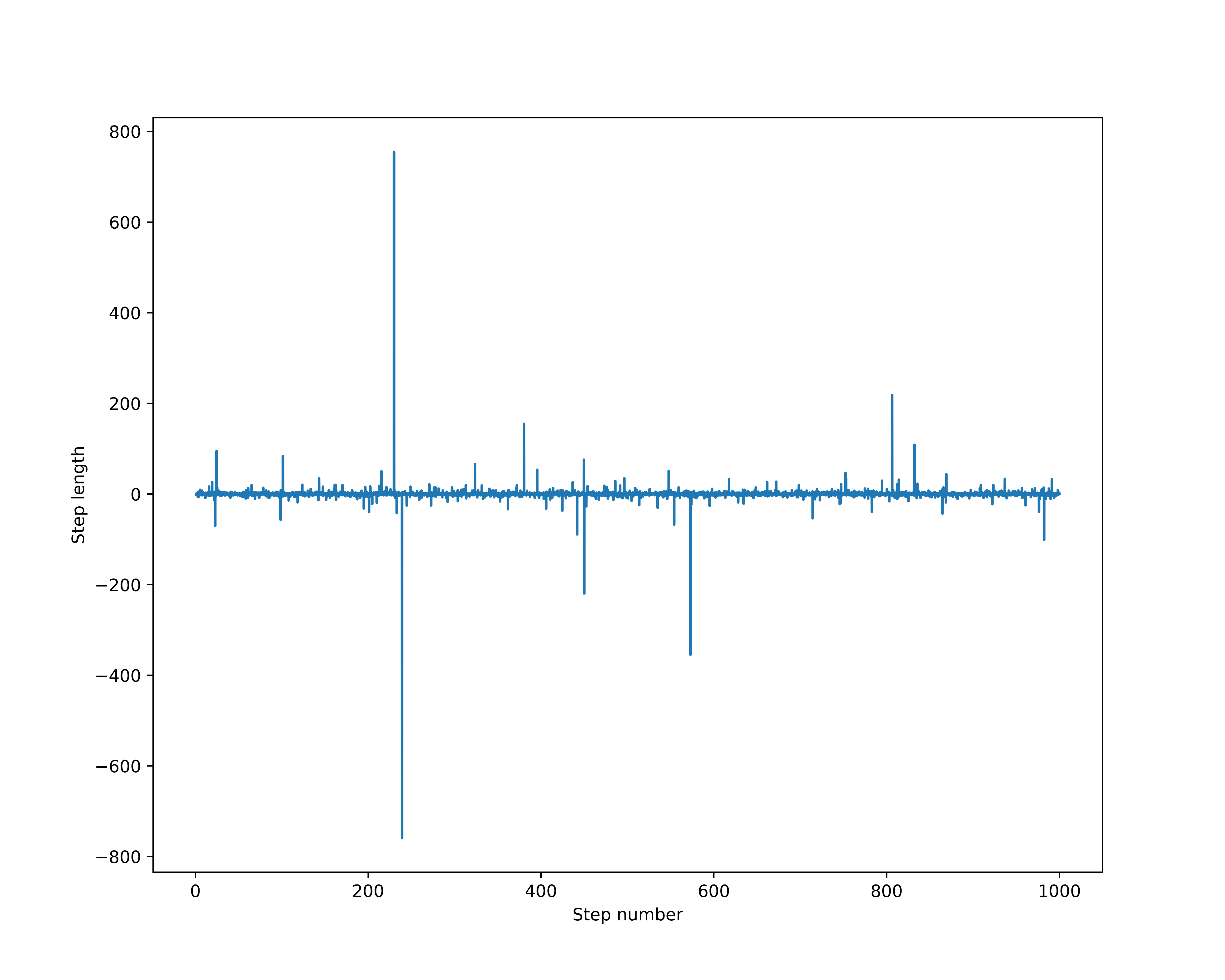}
		\end{minipage}
	}%
	\centering
	\caption{Brownian motion and L\'{e}vy Flight.}\label{fig1}
\end{figure*}

In the second stage, the algorithm transitions into the exploration stage, where juvenile calico salmon migrate from rivers to the ocean. The primary objectives of juvenile calico salmon are to search for food and evade natural predators. The increase in population energy not only changes their feeding behavior but also enhances their ability to evade predators. The mathematical expression for this behavior is represented by Eq. (6). The combination of Brownian motion and L\'{e}vy Flight strategy allows calico salmon to more efficiently pursue prey and evade predators. As shown in Fig. \ref{fig1}, compared to Brownian motion, L\'{e}vy Flight enables the calico salmon to perform longer jumps, while Brownian motion allows them to cover a wider search area. 
\begin{equation}
	\begin{array}{c}
		{x_2}\left( {t + 1} \right) = {x_{best}}\left( t \right) \times {E_1} \\ 
		+ a \times \left( {\left( {UB - LB} \right) \times rand + LB} \right) \\ 
		+ {x_M}\left( t \right) \times {E_1} \times {R_L}\left( \alpha  \right), \\ 
	\end{array}
\end{equation}
where $ x_{2}(t+1) $ represents the solution generated in the next iteration $ t+1 $ by the second stage $ (x_{2}) $. $ a $ is a constant, $ a=1.5 $; and $ R_{L}(\alpha) $ is the L\'{e}vy Flight distribution function, which is calculated using Eq. (7) \cite{bib5}.
\begin{equation}
	{R_L}\left( \alpha  \right) = s \times \frac{m}{{{{\left| n \right|}^{\frac{1}{\alpha }}}}},
\end{equation}
where $ s $ is a constant value fixed to 0.01, $ m $, and $ n $ are two normally distributed quantities defined in Eq. (8), Eq. (9), and the standard deviation $ \sigma_{m} $ and variance $ \sigma_{n} $ defined in Eq. (10).
\begin{equation}
	m = Normal\left( {0,\sigma _m^2} \right),
\end{equation}
\begin{equation}
	n = Normal\left( {0,\sigma _n^2} \right),
\end{equation}
\begin{equation}
	\begin{array}{l}
		{\sigma _x} = {\left[ {\frac{{\Gamma \left( {1 + \alpha } \right) \times \sin \left( {\frac{{\pi \alpha }}{2}} \right)}}{{\Gamma \left( {\frac{{1 + \alpha }}{2}} \right) \times \alpha  \times {2^{\frac{{\left( {\alpha  - 1} \right)}}{2}}}}}} \right]^{\frac{1}{\alpha }}},\:{\sigma _n} = 1,\: \\ 
		and\:\alpha  = 1.5, \\ 
	\end{array}
\end{equation}
$ x_{M}(t) $ denotes the mean of the location of the current solutions connected at $ t^{th} $ iteration, which is calculated by Eq. (11):
\begin{equation}
	{x_M}\left( t \right) = \frac{1}{N}\sum\nolimits_{i = 1}^N {{x_i}\left( t \right)}.  
\end{equation}

\subsubsection{Migrating stage}
In the third stage, when calico salmon reach maturity and the algorithm transitions into the exploitation phase, it becomes essential to prevent being trapped in local optima. To address this, a suitable escape strategy is employed based on population energy, which imitates the behavior of calico salmon escaping from predators. The mathematical formulation of this strategy is presented in Eq. (12).
\begin{equation}
	\begin{array}{c}
		{x_3}\left( {t + 1} \right) = {x_{best}}\left( t \right) \\ 
		+ \left( {{x_{best}}\left( t \right) - {x_R}\left( t \right)} \right) \times {E_1} \times rand \\ 
		+ \left( {{x_R}\left( t \right) - {x_R}\left( {t - 1} \right)} \right){E_2} \times b, \\ 
	\end{array}
\end{equation}
where $ x_{3}(t+1) $ represents the solution generated in the next iteration $ t+1 $ by the third stage $ (x_{3}) $. $ E_{1} $ gradually reach the peak point and then the population energy will shift to the stage $ E_{2} $; $ x_{R}(t-1) $ denotes a random solution taken in the range of $ [ {1,N} ] $ at the $ (t-1)^{th} $ iteration; and $ b $ denotes the control parameter that controls the movement of salmon while avoiding natural predators and is represented by the following formula (13).
\begin{equation}
	b = 1 - \frac{{t}}
	{T}.
\end{equation}

\subsubsection{Mating stage}
In the fourth stage, the remaining calico salmon reproduce their offspring. Although population energy is nearly depleted at this point, they still need to navigate the environment, avoid predators, and locate an appropriate spawning site for fertilization. The mathematical expression for modeling this behavior is provided in Eq. (14).
\begin{equation}
	\begin{array}{c}
		{x_4}\left( {t + 1} \right) = {x_{best}}\left( t \right) \\ 
		+ c\left( {{x_R}\left( t \right) - {x_R}\left( {t - 1} \right)} \right) \times rand \\ 
		+ {E_2} \times {x_R}\left( t \right) \times {R_L}\left( \alpha  \right) \times d, \\ 
	\end{array}
\end{equation}
where $ x_{4}(t+1) $ represents the solution generated in the next iteration $ t+1 $ by the fourth stage $ (x_{4}) $. In formula (15), $ c $ is the control parameter for finding a suitable position during spawning, represented by the following formula (15); $ E_2 $ has been gradually decaying, and the energy of the previous generation of the population tends to be almost zero after the successful fertilization of the calico salmon eggs. $ d $ is the control parameter for the last position update at the end of the spawning, which ranges between -1 and 1.
\begin{equation}
	c = A \times rand - B,
\end{equation}
where $ A=2 $, $ B=1 $.

\begin{figure*}[htbp]%
	\centering
	\includegraphics[width=0.65\textwidth]{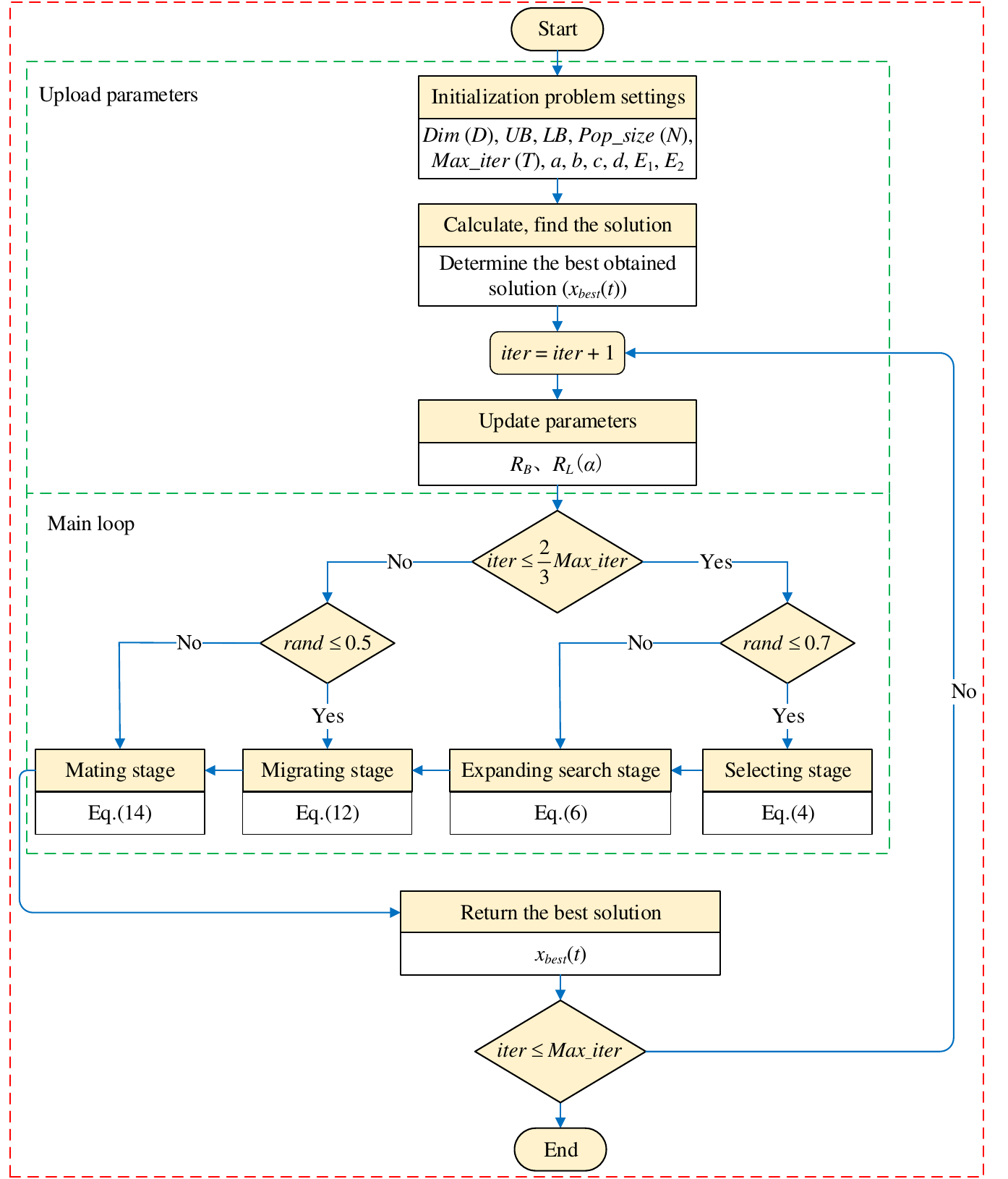}
	\caption{Flowchart of the proposed method.}\label{fig2}
\end{figure*}

In summary, CSMA formulates four different stages for the search strategy: selecting exploration space, expanding the exploration area, reducing the exploration space to focus the search, and ultimately finding the optimal solution. The algorithm iterates until the stopping criteria are met. The flowchart of this method is shown in Fig. \ref{fig2}.

It is noted that the computational complexity of the proposed method is of $ O(CSMA)=O(N*(T*D+1)) $. where, $ T $ number of iterations, $ N $ presents the number of used solutions, and $ D $ presents the solution size.

\section{Experimental results and discussion}
In this section, the newly developed meta-heuristic CSMA is applied to solve 23 CEC-2005 test functions. During the evaluation of CSMA, the numerical experimental results of CSMA are compared with AO \cite{bib6}, WOA \cite{bib7}, HBA \cite{bib8}, ZOA \cite{bib9}, SSA \cite{bib10}, SFO \cite{bib11}, and MSA \cite{bib12}. We use uniform metrics to evaluate the results obtained by CSMA and other competing algorithms. For example, the average of the results (Average) and standard deviation (STD), were run a total of 20 times independently for each test problem. Moreover, a non-parametric statistical test called the Friedman’s mean rank test \cite{bib13} is applied for a fair comparison.

\subsection{Performance comparison through CEC-2005 test functions}
The performance of the CSMA is evaluated using 23 CEC-2005 test functions \cite{bib14}. This set of functions can be divided into three categories. The uni-model functions: These functions have a unique solution, and they are used to test the ability of the algorithm to find a global optimum. Multi-modal functions: These functions have multiple solutions, and they are used to evaluate the algorithm's ability to handle multiple optima. Fixed-dimension functions: These functions are used to test the algorithm's exploration capability in low-dimensional spaces.
\begin{sidewaystable*}
	\caption{The statistical results of the CSMA and other approaches in dealing with the $ (F_1-F_{13}) $}\label{tab1}
	\resizebox{\textwidth}{!}{
	\begin{tabular*}{\textheight}{@{\extracolsep\fill}lcccccccc}
		\toprule%
		Fun No. & \multicolumn{6}{@{}c@{}}{Comparative methods} \\\cmidrule{2-9}%
		Measure & CSMA & WOA & SSA & ZOA & HBA & AO & SFO & MSA  \\
		\midrule
		$ F_1 $ & ~ & ~ & ~ & ~ & ~ & ~ & ~ &   \\  
		Average & 0.0000E+00  & 9.6280E-186 & 6.8171E-11 & 0.0000E+00 & 1.0770E-201 & 4.4349E-07 & 3.0083E-06 & 2.8297E-06  \\ 
		STD & 0.0000E+00 & 0.0000E+00 & 1.5952E-10 & 0.0000E+00 & 0.0000E+00 & 8.3759E-07 & 5.5303E-06 & 2.7533E-06  \\ 
		$ F_3 $ & ~  & ~ & ~ & ~ & ~ & ~ & ~ &   \\ 
		Average & 0.0000E+00  & 2.5550E+04 & 1.4949E-06 & 0.0000E+00 & 3.1233E-129 & 5.7433E-04 & 1.5918E-03 & 3.7903E+03  \\ 
		STD & 0.0000E+00  & 7.8803E+03 & 4.3663E-06 & 0.0000E+00 & 1.3316E-128 & 9.6790E-04 & 3.4021E-03 & 1.9757E+03  \\ 
		$ F_4 $ & ~ & ~ & ~ & ~ & ~ & ~ & ~ &   \\ 
		Average & 7.1661E-169  & 8.3334E+01 & 6.8210E-05 & 1.6663E-244 & 6.9859E-75 & 1.4201E-78 & 1.5028E-04 & 2.7086E+01  \\  
		STD & 0.0000E+00  & 4.2210E+00 & 9.5781E-05 & 0.0000E+00 & 2.0033E-74 & 3.4827E-78 & 1.8382E-04 & 4.8531E+00  \\ 
		$ F_5 $ & ~ & ~ & ~ & ~ & ~ & ~ & ~ &   \\ 
		Average & 9.6690E+00  & 2.7542E+01 & 9.4681E-08 & 2.7958E+01 & 2.7569E+01 & 4.8066E-04 & 3.1514E-05 & 1.6893E+02  \\ 
		STD & 1.3153E+01  & 6.7776E-01 & 1.9184E-07 & 7.3717E-01 & 1.4294E+00 & 1.3088E-03 & 3.9916E-05 & 3.1810E+02  \\ 
		$ F_6 $ & ~ & ~ & ~ & ~ & ~ & ~ & ~ &   \\ 
		Average & 6.9651E-05  & 1.4902E-01 & 2.8632E-10 & 2.1333E+00 & 2.4399E+00 & 2.0593E-06 & 2.8188E+00 & 2.8720E-06  \\  
		STD & 6.5870E-05  & 1.2870E-01 & 6.1122E-10 & 6.0608E-01 & 5.6234E-01 & 4.1230E-06 & 1.4118E+00 & 4.4416E-06  \\ 
		$ F_7 $ & ~ & ~ & ~ & ~ & ~ & ~ & ~ &   \\ 
		Average & 5.8707E-05  & 1.7171E-03 & 1.7387E-04 & 4.0574E-05 & 2.5535E-04 & 2.7455E-04 & 2.0581E-03 & 1.3947E-01  \\ 
		STD & 4.7874E-05  & 2.1474E-03 & 1.6693E-04 & 2.7729E-05 & 1.8383E-04 & 2.5577E-04 & 1.3269E-03 & 5.4616E-02  \\ 
		$ F_8 $ & ~ & ~  & ~ & ~ & ~ & ~ & ~ &   \\ 
		Average & -1.2569E+04  & -8.6062E+03 & -6.6883E+03 & -6.7792E+03 & -7.5580E+03 & -1.2399E+04 & -3.9610E+03 & -9.4283E+03  \\ 
		STD & 5.3268E-04  & 3.5303E+02 & 1.4197E+03 & 6.4466E+02 & 8.0830E+02 & 4.4044E+02 & 5.8025E+02 & 5.7681E+02  \\
		$ F_9 $ & ~ & ~  & ~ & ~ & ~ & ~ & ~ &   \\  
		Average & 0.0000E+00 & 4.2773E+01 & 9.5449E-09 & 0.0000E+00 & 0.0000E+00 & 2.3978E-06 & 2.9387E-04 & 1.2855E+01  \\  
		STD & 0.0000E+00 & 7.5418E+01 & 2.9730E-08 & 0.0000E+00 & 0.0000E+00 & 3.2015E-06 & 6.9950E-04 & 3.9241E+00  \\ 
		$ F_{10} $ & ~  & ~ & ~ & ~ & ~ & ~ & ~ &   \\ 
		Average & 4.4409E-16  & 3.9968E-15 & 3.9073E-06 & 1.8652E-15 & 3.9968E-15 & 3.9968E-15 & 9.6098E-04 & 1.7168E-03  \\ 
		STD & 0.0000E+00  & 3.5527E-15 & 6.3661E-06 & 1.7405E-15 & 0.0000E+00 & 0.0000E+00 & 1.2108E-03 & 4.4569E-04  \\ 
		$ F_{11} $ & ~ & ~  & ~ & ~ & ~ & ~ & ~ &   \\ 
		Average & 0.0000E+00  & 1.7446E-03 & 3.8740E-10 & 0.0000E+00 & 0.0000E+00 & 0.0000E+00 & 2.6657E-07 & 3.3523E+00  \\  
		STD & 0.0000E+00  & 7.6046E-03 & 1.6613E-09 & 0.0000E+00 & 0.0000E+00 & 0.0000E+00 & 8.2694E-07 & 2.3103E+00  \\ 
		$ F_{12} $ & ~ & ~  & ~ & ~ & ~ & ~ & ~ &   \\  
		Average & 1.7068E-06  & 1.1861E+00 & 1.7841E-05 & 1.2206E-01 & 1.3232E-01 & 6.2568E-03 & 4.8863E-01 & 5.7073E-02  \\  
		STD & 1.6973E-06  & 3.0541E+00 & 1.3801E-05 & 8.2130E-02 & 1.1545E-01 & 2.3899E-04 & 4.0573E-01 & 1.2049E-01  \\ 
		$ F_{13} $ & ~ & ~  & ~ & ~ & ~ & ~ & ~ &   \\  
		Average & 1.6037E-06  & 6.4255E-01 & 4.9251E-09 & 1.8506E+00 & 1.7130E+00 & 1.7298E-06 & 1.6178E-07 & 6.7908E-07  \\  
		STD & 1.9695E-06  & 2.9440E-01 & 9.2094E-09 & 3.6413E-01 & 1.4441E-01 & 2.6322E-06 & 2.0881E-07 & 1.1733E-06 \\ 
		Mean & 2.5300E+00  & 5.7800E+00 & 4.3300E+00 & 4.3000E+00 & 5.1000E+00 & 4.1400E+00 & 6.7700E+00 & 6.69000E+00  \\  
		\bottomrule
	\end{tabular*}}
\end{sidewaystable*}
\begin{sidewaystable*}
	\caption{The statistical results of the CSMA and other approaches in dealing with the $ (F_{14}-F_{20}) $}\label{tab2}
	\resizebox{\textwidth}{!}{
	\begin{tabular*}{\textheight}{@{\extracolsep\fill}lcccccccc}
		\toprule%
		Fun No. & \multicolumn{8}{@{}c@{}}{Comparative methods} \\\cmidrule{2-9}%
		Measure & CSMA  & WOA & SSA & ZOA & HBA & AO & SFO & MSA  \\
		\midrule
		$ F_{14} $ & ~ & ~ & ~ & ~ & ~ & ~ & ~ &   \\ 
		Average & 9.9800E-01 & 1.2960E+00 & 9.5686E+00 & 2.8708E+00 & 1.0972E+00 & 3.0732E+00 & 4.1215E+00 & 3.8838E+00  \\ 
		STD & 1.7378E-16 & 5.5285E-01 & 4.6092E+00 & 2.4877E+00 & 4.3242E-01 & 2.9719E+00 & 2.5798E+00 & 4.1546E+00  \\ 
		$ F_{15} $ & ~ & ~ & ~ & ~ & ~ & ~ & ~ &   \\  
		Average & 3.0749E-04 & 6.6353E-04 & 3.2346E-04 & 3.1694E-04 & 6.8361E-04 & 1.1363E-03 & 5.7076E-04 & 2.7453E-03  \\  
		STD & 2.7428E-19 & 5.4598E-04 & 6.9414E-05 & 3.2149E-05 & 5.2394E-04 & 4.8144E-04 & 1.8183E-04 & 5.8784E-03  \\ 
		$ F_{16}$ & ~ & ~ & ~ & ~ & ~ & ~ & ~ &   \\  
		Average & -1.0316E+00 & -1.0316E+00 & -9.9082E-01 & -1.0316E+00 & -1.0316E+00 & -1.0272E+00 & -1.0062E+00 & -1.0316E+00  \\ 
		STD & 0.0000E+00 & 0.0000E+00 & 1.7788E-01 & 5.6606E-11 & 0.0000E+00 & 6.1887E-03 & 1.6612E-02 & 7.0012E-12  \\ 
		$ F_{17} $ & ~ & ~ & ~ & ~ & ~ & ~ & ~ &   \\  
		Average & 3.9789E-01 & 3.9789E-01 & 3.9789E-01 & 3.9789E-01 & 3.9789E-01 & 4.3308E-01 & 4.2414E-01 & 3.9789E-01  \\  
		STD & 0.0000E+00 & 2.0876E-14 & 8.8896E-15 & 1.3375E-08 & 0.0000E+00 & 3.2830E-02 & 3.6509E-02 & 1.0226E-12  \\ 
		$ F_{18} $ & ~ & ~ & ~ & ~ & ~ & ~ & ~ &   \\  
		Average & 3.0000E+00 & 3.0000E+00 & 7.0500E+00 & 5.7000E+00 & 9.7500E+00 & 4.7101E+00 & 4.2904E+00 & 3.0000E+00  \\  
		STD & 2.5278E-15 & 5.4217E-15 & 9.6409E+00 & 8.1000E+00 & 1.8852E+01 & 2.9279E+00 & 1.3236E+00 & 8.9868E-10  \\ 
		$ F_{19} $ & ~ & ~ & ~ & ~ & ~ & ~ & ~ &   \\  
		Average & -3.0524E+00 & -3.0524E+00 & -3.0524E+00 & -3.0524E+00 & -3.0524E+00 & -2.3137E+00 & -2.9832E+00 & -3.0524E+00  \\ 
		STD & 0.0000E+00 & 1.9860E-15 & 4.3460E-14 & 1.1147E-06 & 0.0000E+00 & 4.5345E-01 & 6.1568E-02 & 4.7559E-08  \\ 
		$ F_{20} $ & ~ & ~ & ~ & ~ & ~ & ~ & ~ &   \\  
		Average & -3.2804E+00 & -3.2712E+00 & -3.2745E+00 & -3.3099E+00 & -3.2863E+00 & -2.2975E+00 & -2.7173E+00 & -3.3101E+00  \\  
		STD & 5.6721E-02 & 6.3754E-02 & 5.8258E-02 & 3.6253E-02 & 5.4495E-02 & 5.8601E-01 & 2.4655E-01 & 3.5676E-02  \\
		$ F_{21} $ & ~ & ~ & ~ & ~ & ~ & ~ & ~ &   \\ 
		Average & -1.0153E+01 & -8.2477E+00 & -1.0153E+01 & -9.8983E+00 & -9.0248E+00 & -1.0147E+01 & -1.0153E+01 & -6.3895E+00  \\  
		STD & 0.0000E+00 & 2.6449E+00 & 3.2657E-07 & 1.1111E+00 & 2.6862E+00 & 3.8656E-03 & 6.0383E-05 & 3.4970E+00  \\ 
		$ F_{22} $ & ~ & ~ & ~ & ~ & ~ & ~ & ~ &   \\  
		Average & -1.0403E+01 & -7.0780E+00 & -9.8714E+00 & -9.3395E+00 & -8.8755E+00 & -1.0399E+01 & -1.0403E+01 & -4.5483E+00  \\  
		STD & 0.0000E+00 & 3.4064E+00 & 1.5946E+00 & 2.1259E+00 & 3.0548E+00 & 3.6946E-03 & 7.8929E-05 & 2.9825E+00  \\ 
		$ F_{23} $ & ~ & ~ & ~ & ~ & ~ & ~ & ~ &   \\ 
		Average & -1.0536E+01 & -7.0225E+00 & -1.0536E+01 & -9.9956E+00 & -8.8612E+00 & -1.0530E+01 & -1.0536E+01 & -6.1379E+00  \\  
		STD & 0.0000E+00 & 3.9317E+00 & 5.35600E-05 & 1.6224E+00 & 2.9016E+00 & 7.24220E-03 & 1.4802E-04 & 3.67770E+00  \\ 
		Mean & 2.3100E+00 & 3.8800E+00 & 3.7300E+00 & 4.9100E+00 & 3.3200E+00 & 8.1000E+00 & 7.6100E+00 & 5.85000E+00  \\  
		\bottomrule
	\end{tabular*}}
\end{sidewaystable*}

\begin{table*}[htbp]
	\caption{\centering{$ p-values $ of the  Friedman’s mean rank test over 20 runs $ (F_{1}-F_{7}) $}}\label{tab3}
	\resizebox{\textwidth}{!}{
		\begin{tabular*}{\textheight}{@{\extracolsep\fill}lccccccccc}
			\toprule
			Result & \multicolumn{9}{@{}c@{}}{Comparative methods} \\\cmidrule{2-10}
			Measure & CSMA  & WOA & SSA & ZOA & HBA & AO & SFO & MSA  & $ p-values $ \\
			\midrule
			Mean & 2.28 & 3.59 & 3.47 & 4.66 & 3.14 & 6.96 & 6.55 & 5.36 & \multirow{2}{*}{$ < 0.001 $}\\
			Ranking & 1 & 4 & 3 & 5 & 2 & 8 & 7 & 6 & \\
			\bottomrule
	\end{tabular*}}
\end{table*}
\begin{table*}[htbp]
	\caption{\centering{$ p-values $ of the  Friedman’s mean rank test over 20 runs $ (F_{8}-F_{13}) $}}\label{tab4}
	\resizebox{\textwidth}{!}{
		\begin{tabular*}{\textheight}{@{\extracolsep\fill}lccccccccc}
			\toprule
			Result & \multicolumn{9}{@{}c@{}}{Comparative methods} \\\cmidrule{2-10}
			Measure & CSMA  & WOA & SSA & ZOA & HBA & AO & SFO & MSA  & $ p-values $ \\
			\midrule
			Mean & 2.59 & 6.07 & 3.29 & 3.16 & 4.57 & 4.02 & 5.64 & 6.679 & \multirow{2}{*}{$ < 0.001 $}\\
			Ranking & 1 & 7 & 3 & 2 & 5 & 4 & 6 & 8 & \\
			\bottomrule
	\end{tabular*}}
\end{table*}
\begin{table*}[htbp]
	\caption{\centering{$ p-values $ of the  Friedman’s mean rank test over 20 runs $ (F_{14}-F_{23}) $}}\label{tab5}
	\resizebox{\textwidth}{!}{
		\begin{tabular*}{\textheight}{@{\extracolsep\fill}lccccccccc}
			\toprule
			Result & \multicolumn{9}{@{}c@{}}{Comparative methods} \\\cmidrule{2-10}
			Measure & CSMA  & WOA & SSA & ZOA & HBA & AO & SFO & MSA  & $ p-values $ \\
			\midrule
			Mean & 2.28 & 3.59 & 3.47 & 4.66 & 3.14 & 6.96 & 6.55 & 5.36 & \multirow{2}{*}{$ < 0.001 $}\\
			Ranking & 1 & 4 & 3 & 5 & 2 & 8 & 7 & 6 & \\
			\bottomrule
	\end{tabular*}}
\end{table*}

\subsection{Comparison of CSMA with some latest algorithms}
Table \ref{tab1} and Table \ref{tab2} provide the statistical results of CSMA compared to other competing algorithms (Average and STD values). From the results of the numerical experiments, it can be observed that CSMA exhibits high efficiency in successfully finding global optima for both uni-model, multi-modal, and fixed-dimension functions. $ {F_1} - {F_7} $ primarily demonstrate a fast convergence rate. For $ {F_8} - {F_{13}} $, the competitive algorithms in the same category perform well both in convergence speed and convergence accuracy CSMA. CSMA wins 5 times and the other competitive algorithms are ZOA and SSA respectively. The performance of CSMA will be more obvious in $ F_{14} $, $ F_{19} $, $ F_{20} $, $ F_{21} $, $ F_{22} $, and in $ F_{16} $, $ F_{18} $, $ F_{23} $, CSMA is on par with other competing algorithms, especially MSA and SFO, which are highly competitive. For $ F_{15} $, $ F_{17} $, CSMA has a slower decline in the early stage, but it can catch up with other competing algorithms in the middle and late iterations until the optimal solution is found. 

\subsection{Friedman’s mean rank test results}
This section supplements the Friedman’s mean rank test to confirm the significance of the experimental results. The $ p-values $ in Table \ref{tab3}, Table \ref{tab4}, and Table \ref{tab5} represent the results of the Friedman’s mean rank test. The experimental results for both uni-modal functions, multi-modal functions, and fixed-dimensional test functions all indicate that CSMA numerical experiments have statistical significance. The null hypothesis is rejected, further demonstrating that CSMA has a distinct advantage over other competitive algorithms.

\section{Conclusion}
In this paper, a new meta-heuristic technique, Calico Salmon Migration Algorithm (CSMA) is proposed, which simulates the natural behavior of calico salmon moving rivers into the sea and migrating for spawning in nature. In CSMA, the optimization process is represented by four stages: selecting stage, expanding search stage, migrating phase, and mating stage. CSMA combines L\'{e}vy flight, Brownian motion, evolutionary operators, and heuristic algorithms to enhance its performance. The statistical results from numerical experiments indicate that CSMA either outperforms other well-known meta-heuristic algorithms or is at least comparable to them.

\section*{Acknowledgment}
This work was supported by the National Natural Science Foundation of China (No: 11901484).



%

\end{document}